# On the Hierarchical Preconditioning of the Combined Field Integral Equation

Simon B. Adrian[1,2], Francesco P. Andriulli[2], and Thomas F. Eibert[1]
[1]Lehrstuhl für Hochfrequenztechnik, Technische Universität München, Munich, Germany
[2]Institut Mines-Télécom / Télécom Bretagne, Technopôle Brest-Iroise, Brest, France

*Abstract*—**This paper analyzes how hierarchical bases preconditioners constructed for the Electric Field Integral Equation (EFIE) can be effectively applied to the Combined Field Integral Equation (CFIE). For the case where no hierarchical solenoidal basis is available (e.g., on unstructured meshes), a new scheme is proposed: the CFIE is implicitly preconditioned on the solenoidal Helmholtz subspace by using a Helmholtz projector, while a hierarchical non-solenoidal basis is used for the non-solenoidal Helmholtz subspace. This results in a well-conditioned system. Numerical results corroborate the presented theory.**

## I. Introduction

Scattering and radiation problems in electromagnetics are commonly solved by leveraging the Electric Field Integral Equation (EFIE) in the framework of boundary element methods. The EFIE, however, is not well-conditioned: whenever the frequency decreases or the discretization density increases, the condition number grows [1], [2]; these effects are often called low-frequency and dense-discretization breakdown. In addition, the EFIE is not uniquely solvable at resonance frequencies [3].

The low-frequency breakdown has been solved in the past by using a quasi-Helmholtz decomposition of the current, such as the loop-star/loop-tree decomposition, allowing for a separate rescaling of the vector and the scalar potential in frequency [2]. The dense-discretization breakdown was remedied by the introduction of hierarchical basis preconditioners [4]–[7]. Since hierarchical basis preconditioners are usually based on a quasi-Helmholtz decomposition, they cure the low-frequency breakdown as well. The interior resonances have been removed by solving a linear combination of the EFIE and of the Magnetic Field Integral Equation (MFIE): the Combined Field Integral Equation (CFIE) [3]. Since, however, the CFIE is a linear combination of the EFIE and of the MFIE, it is not free from the low-frequency and the dense-discretization breakdown. For this reason, efforts have been made to extend the hierarchical basis preconditioners to the CFIE [8].

The conference contribution [8] reported on the fact that the application of a hierarchical loop (h-loop)/hierarchical non-solenoidal (h-n-sol) basis preconditioner to the CFIE resulted in a well-conditioned equation. Given that the h-loop preconditioner is not capable of preconditioning the EFIE [7], [9], we found the result surprising and decided to perform a theoretical investigation of the problem. This lead to the current letter whose novelty content is twofold: (i) we show on a theoretically sound basis that a direct application of the hierarchical basis preconditioner to the CFIE is possible when the h-loop [4] or the 3-point hierarchical loop (3ph-loop) functions [7] are used as solenoidal basis. (ii) We propose a new hierarchical left preconditioning scheme that differently from any other hierarchical approach to the preconditioning of CFIE (including those communicated in [8] and detailed in [10]) can provide optimal preconditioning *both* on structured and unstructured meshes. Numerical results will show the effectiveness of the new scheme for canonical and realistic examples. Very preliminary results have appeared as conference contribution [11].

## II. Notation and Background

We consider an electromagnetic wave $(\boldsymbol{E}^{\mathrm{i}}, \boldsymbol{H}^{\mathrm{i}})$ impinging on a Perfectly Electrically Conducting (PEC) scatterer $\Omega$ with boundary $\Gamma = \partial\Omega$, which is embedded in a homogeneous medium with permittivity $\varepsilon$ and permeability $\mu$. The impinging wave excites an electric surface current density $\boldsymbol{j}$ generating the scattered wave $(\boldsymbol{E}^{\mathrm{s}}, \boldsymbol{H}^{\mathrm{s}})$. Thus the total field is $(\boldsymbol{E}, \boldsymbol{H}) = (\boldsymbol{E}^{\mathrm{i}} + \boldsymbol{E}^{\mathrm{s}}, \boldsymbol{H}^{\mathrm{i}} + \boldsymbol{H}^{\mathrm{s}})$ and on $\Gamma$ the tangential component of the electric field satisfies the boundary condition $\hat{\boldsymbol{n}} \times \boldsymbol{E} = \boldsymbol{0}$.

The surface current density $\boldsymbol{j}$ is recovered from the tangential component of the incident electric field by solving the EFIE $\mathcal{T}\boldsymbol{j} = -\hat{\boldsymbol{n}} \times \boldsymbol{E}^{\mathrm{i}}$ and from the tangential component of the magnetic field by solving the MFIE $\mathcal{M}\boldsymbol{j} = (\mathcal{I}/2 + \mathcal{K})\,\boldsymbol{j} = \eta\hat{\boldsymbol{n}} \times \boldsymbol{H}^{\mathrm{i}}$ on the surface $\Gamma$ with wave impedance $\eta = \sqrt{\mu/\varepsilon}$, where we assumed that the current is normalized such that $\boldsymbol{j} = \eta\hat{\boldsymbol{n}} \times \boldsymbol{H}$; the operator $\mathcal{T} = \mathcal{T}_{\mathrm{A}} + \mathcal{T}_{\Phi}$ is the EFIE operator with the vector potential $\mathcal{T}_{\mathrm{A}}\boldsymbol{j} = \hat{\boldsymbol{n}} \times \mathrm{i}k \int_\Gamma \frac{\mathrm{e}^{\mathrm{i}k|\boldsymbol{r}-\boldsymbol{r}'|}}{4\pi|\boldsymbol{r}-\boldsymbol{r}'|}\boldsymbol{j}(\boldsymbol{r}')\mathrm{d}\boldsymbol{r}'$, scalar potential $\mathcal{T}_{\Phi}\boldsymbol{j} = -\hat{\boldsymbol{n}} \times \frac{1}{\mathrm{i}k}\nabla_\Gamma \int_\Gamma \frac{\mathrm{e}^{\mathrm{i}k|\boldsymbol{r}-\boldsymbol{r}'|}}{4\pi|\boldsymbol{r}-\boldsymbol{r}'|}\nabla_\Gamma \cdot \boldsymbol{j}(\boldsymbol{r}')\mathrm{d}\boldsymbol{r}'$, and $\mathcal{M} = \mathcal{I}/2 + \mathcal{K}$ is the MFIE operator with $\mathcal{K}\boldsymbol{j} = -\hat{\boldsymbol{n}} \times \left(\nabla_\Gamma \times \int_\Gamma \frac{\mathrm{e}^{\mathrm{i}k|\boldsymbol{r}-\boldsymbol{r}'|}}{4\pi|\boldsymbol{r}-\boldsymbol{r}'|}\boldsymbol{j}(\boldsymbol{r}')\mathrm{d}\boldsymbol{r}'\right)$, where $k$ is the wavenumber and $\hat{\boldsymbol{n}}$ is the exterior surface normal vector. The CFIE is given by $\mathcal{C}\boldsymbol{j} = -\alpha\hat{\boldsymbol{n}} \times \hat{\boldsymbol{n}} \times \boldsymbol{E}^{\mathrm{i}} + \eta(1-\alpha)\hat{\boldsymbol{n}} \times \boldsymbol{H}^{\mathrm{i}}$, where $\mathcal{C} = \alpha\hat{\boldsymbol{n}} \times \mathcal{T} + (1-\alpha)\mathcal{M}$ is the CFIE operator with $\alpha \in ]0, 1[$.

Following a Galerkin approach, the CFIE is commonly discretized using Rao-Wilton-Glisson (RWG) functions $\boldsymbol{f}_i \in X_{\mathrm{RWG}}$ as source and testing functions. Hence, the CFIE system reads

$$\boldsymbol{Z}_{\mathrm{C}}\boldsymbol{i} = \boldsymbol{v}_{\mathrm{C}} \tag{1}$$

where $\left[\boldsymbol{Z}_{\mathrm{C}}\right]_{ij} = \alpha\boldsymbol{Z}_{\mathrm{E}} + (1-\alpha)\eta\boldsymbol{Z}_{\mathrm{M}}$, $\boldsymbol{j} \approx \sum_{i=1}^{N}\left[\boldsymbol{i}\right]_j \boldsymbol{f}_j$, $\boldsymbol{v}_{\mathrm{C}} = \alpha\boldsymbol{v}_{\mathrm{E}} + (1-\alpha)\eta\boldsymbol{v}_{\mathrm{M}}$ with $\boldsymbol{Z}_{\mathrm{E}} = \boldsymbol{Z}_{\mathrm{A}} + \boldsymbol{Z}_{\Phi}$, $\left[\boldsymbol{v}_{\mathrm{E}}\right]_i = -\left(\boldsymbol{f}_i, \hat{\boldsymbol{n}} \times \hat{\boldsymbol{n}} \times \boldsymbol{E}^{\mathrm{i}}\right)_{L^2}$ and $\left[\boldsymbol{Z}_{\mathrm{A}}\right]_{ij} = \left(\boldsymbol{f}_i, \hat{\boldsymbol{n}} \times \mathcal{T}_{\mathrm{A}}\boldsymbol{f}_j\right)_{L^2}$, $\left[\boldsymbol{Z}_{\Phi}\right]_{ij} = \left(\boldsymbol{f}_i, \hat{\boldsymbol{n}} \times \mathcal{T}_{\Phi}\boldsymbol{f}_j\right)_{L^2}$ and $\left[\boldsymbol{Z}_{\mathrm{M}}\right]_{ij} = \left(\boldsymbol{f}_i, \mathcal{M}\boldsymbol{f}_j\right)_{L^2}$, $\left[\boldsymbol{v}_{\mathrm{M}}\right]_i = \left(\boldsymbol{f}_i, \hat{\boldsymbol{n}} \times \boldsymbol{H}^{\mathrm{i}}\right)_{L^2}$.

It has been pointed out that testing the MFIE operator with standard RWG functions is not consistent with the correct discretization framework of the operator and that one way of obtaining a correct and conforming discretization would be to use Buffa-Christiansen (BC) functions $\widetilde{\boldsymbol{f}}_i$ as testing functions [12]. In this way, the conformingly discretized MFIE system reads $\left[\widetilde{\boldsymbol{Z}}_\mathrm{M}\right]_{ij} = \left(\hat{\boldsymbol{n}} \times \widetilde{\boldsymbol{f}}_i, \mathcal{M}\boldsymbol{f}_j\right)_{L^2}$ and $\left[\boldsymbol{v}_\mathrm{M}\right]_i = \left(\hat{\boldsymbol{n}} \times \widetilde{\boldsymbol{f}}_i, \hat{\boldsymbol{n}} \times \boldsymbol{H}^\mathrm{i}\right)_{L^2}$. Accordingly, the modified CFIE matrix reads

$$\left[\widetilde{\boldsymbol{Z}}_\mathrm{C}\right]_{ij} = \alpha \boldsymbol{Z}_\mathrm{E} + (1-\alpha)\eta \boldsymbol{G}_\mathrm{ff}\boldsymbol{G}^{-1}_{\mathrm{n}\times\tilde{\mathrm{f}},\mathrm{f}}\boldsymbol{Z}_\mathrm{M}, \tag{2}$$

and $\widetilde{\boldsymbol{v}}_\mathrm{C} = \alpha \boldsymbol{v}_\mathrm{E} + (1-\alpha)\eta \boldsymbol{G}_\mathrm{ff}\boldsymbol{G}^{-1}_{\mathrm{n}\times\tilde{\mathrm{f}},\mathrm{f}}\boldsymbol{v}_\mathrm{M}$ with the Gram matrices $\left[\boldsymbol{G}_\mathrm{ff}\right]_{ij} = \left(\boldsymbol{f}_i, \boldsymbol{f}_j\right)_{L^2}$ and $\left[\boldsymbol{G}_{\mathrm{n}\times\tilde{\mathrm{f}},\mathrm{f}}\right]_{ij} = \left(\hat{\boldsymbol{n}} \times \widetilde{\boldsymbol{f}}_i, \boldsymbol{f}_j\right)_{L^2}$.

It is well known that the EFIE is ill-conditioned with $\operatorname{cond}(\boldsymbol{Z}_\mathrm{E}) = O\left(1/(hk)^2\right)$, where $h$ is the average edge length of the mesh. Hierarchical basis preconditioners that cure the EFIE ill-conditioning consist of two sets of basis functions, hierarchical solenoidal (h-sol) and h-n-sol functions forming a quasi-Helmholtz decomposition of $X_\mathrm{RWG}$.

Before we study the hierarchical preconditioners, we consider the loop-star quasi-Helmholtz decomposition of the EFIE. Let $\boldsymbol{\Lambda}_i$ be the loop functions and $\boldsymbol{\Sigma}_i$ be the star functions. Then we have $\boldsymbol{j} \approx \sum_{i=1}^{N}\left[\boldsymbol{i}\right]_j \boldsymbol{f}_j = \sum_{i=1}^{N_\Lambda}\left[\boldsymbol{i}_\Lambda\right]_j \boldsymbol{\Lambda}_j + \sum_{i=1}^{N_\Sigma}\left[\boldsymbol{i}_\Sigma\right]_j \boldsymbol{\Sigma}_j$. Since the loop and star functions are defined as linear combinations of RWG functions, we can define transformation matrices $\boldsymbol{\Lambda}$ and $\boldsymbol{\Sigma}$ such that $\boldsymbol{i} = \boldsymbol{\Lambda}\boldsymbol{i}_\Lambda + \boldsymbol{\Sigma}\boldsymbol{i}_\Sigma$. We define $\boldsymbol{T} = \left[\boldsymbol{\Lambda}/\sqrt{k} \quad \boldsymbol{\Sigma}\sqrt{k}\right]$ and the EFIE system matrix in the loop-star basis is given by $\boldsymbol{Z}_{\mathrm{E,qH}} := \boldsymbol{T}^\mathsf{T}\boldsymbol{Z}_\mathrm{E}\boldsymbol{T} = \begin{bmatrix}\boldsymbol{Z}_{\Lambda\Lambda} & \boldsymbol{Z}_{\Lambda\Sigma}\\ \boldsymbol{Z}_{\Sigma\Lambda} & \boldsymbol{Z}_{\Sigma\Sigma}\end{bmatrix}$. In the low-frequency limit, $\boldsymbol{Z}_{\mathrm{E,qH}}$ becomes block diagonal, that is, $\boldsymbol{Z}_{\Lambda\Sigma}$ and $\boldsymbol{Z}_{\Sigma\Lambda}$ vanish; thus, $\boldsymbol{Z}_{\Lambda\Lambda}$ and $\boldsymbol{Z}_{\Sigma\Sigma}$ can be identified as the sources of the ill-conditioning. In fact, it has been shown that the matrices $\boldsymbol{Z}_{\Lambda\Lambda,0} = \lim_{k\to 0}\boldsymbol{Z}_{\Lambda\Lambda}$ and $\boldsymbol{Z}_{\Sigma\Sigma,0} = \lim_{k\to 0}\boldsymbol{Z}_{\Sigma\Sigma}$ can be expressed as $\left[\boldsymbol{Z}_{\Lambda\Lambda,0}\right]_{ij} = \left(\lambda_i, \mathcal{W}\lambda_j\right)_{L_2}$ and $\left[\boldsymbol{Z}_{\Sigma\Sigma,0}\right]_{ij} = \left(\sigma_i, \mathcal{V}\sigma_j\right)_{L_2}$, where $\lambda_i$ are the nodal functions (i.e., $\lambda_i(\boldsymbol{r}) = 1$ for $\boldsymbol{r} \in n_i$, $\lambda_i(\boldsymbol{r}) = 0$ for $\boldsymbol{r} \in n_j \neq n_i$, where $n_i \subset \Gamma$ is the $i$th node of the mesh), $\sigma_i = \operatorname{div}\boldsymbol{\Sigma}_i$, and $\boldsymbol{\Lambda}_i = \nabla \times \hat{\boldsymbol{n}}\lambda_i$. Moreover the operators are defined as $\mathcal{W}f = \hat{\boldsymbol{n}}_{\boldsymbol{r}} \cdot \nabla \times \int_\Gamma \frac{1}{4\pi|\boldsymbol{r}-\boldsymbol{r}'|}\nabla' \times \hat{\boldsymbol{n}}_{\boldsymbol{r}'} f(\boldsymbol{r}')\mathrm{d}S(\boldsymbol{r}')$ and $\mathcal{V}f = \int_\Gamma \frac{1}{4\pi|\boldsymbol{r}-\boldsymbol{r}'|}f(\boldsymbol{r}')\mathrm{d}S(\boldsymbol{r}')$ [7].

## III. Spectral Analysis and Preconditioning

The operators $\mathcal{W}$ and $\mathcal{V}$ are known to suffer from the dense-discretization breakdown [7]. If we precondition these operators successfully, we can precondition the EFIE. Different hierarchical bases have been introduced in the past. For preconditioning $\mathcal{T}_\Phi$, the bases developed are all equivalent, and we choose as a representative the basis developed in [5]. For $\mathcal{T}_\mathrm{A}$, two hierarchical bases have been presented: the h-loops $\boldsymbol{\Lambda}_{\mathrm{H}i}$ [4] and the 3ph-loops $\boldsymbol{\Lambda}_{\mathrm{T}i}$ [7]. Only the latter can cure, however, the dense-discretization breakdown of the EFIE [7]. Both $\boldsymbol{\Lambda}_{\mathrm{H}i}$ and $\boldsymbol{\Lambda}_{\mathrm{T}i}$ require a structured mesh: a mesh that is obtained by iteratively refining, for example dyadically, $J$ times a coarse mesh. Accordingly, $\boldsymbol{\Lambda}_{\mathrm{H}i}$ and $\boldsymbol{\Lambda}_{\mathrm{T}i}$ are defined on $J$ hierarchical levels (for details, see [7]). These functions are related to scalar functions by $\boldsymbol{\Lambda}_{\mathrm{H}i} = \nabla \times \hat{\boldsymbol{n}}\lambda_{\mathrm{H}i}$ and $\boldsymbol{\Lambda}_{\mathrm{T}i} = \nabla \times \hat{\boldsymbol{n}}\lambda_{\mathrm{T}i}$, where $\lambda_{\mathrm{H}i}$ are the hierarchical nodal (h-nodal) functions presented by Yserentant [13] and $\lambda_{\mathrm{T}i}$ are the 3-point hierarchical nodal (3ph-nodal) functions by Stevenson [14]. As we did for the loop basis, we can define transformation matrices $\boldsymbol{\Lambda}_\mathrm{H}$ and $\boldsymbol{\Lambda}_\mathrm{T}$ for the h-loop and 3ph-loop functions, respectively.

The operator $\mathcal{W}$ maps from the Sobolev space $H^{1/2}(\Gamma)$ to $H^{-1/2}(\Gamma)$ and induces an inner product on $H^{1/2}/\mathbb{R}$ (for a definition of the Sobolev spaces, see [15]). The 3ph-nodal functions $\lambda_{\mathrm{T}i}$ are $H^{1/2}$-stable when rescaled appropriately, and therefore, the 3ph-loop functions can precondition the vector potential [7] resulting in $\operatorname{cond}\left(\boldsymbol{D}^1_\Lambda\boldsymbol{\Lambda}^\mathsf{T}_\mathrm{T}\boldsymbol{Z}_\mathrm{A}\boldsymbol{\Lambda}_\mathrm{T}\boldsymbol{D}^1_\Lambda\right) \lesssim 1$, where $\boldsymbol{D}_{\boldsymbol{\Lambda}_\mathrm{T}}$ is a diagonal matrix with entries $\left[\boldsymbol{D}^s_\Lambda\right]_{ii} = 2^{sl(i)}$, where $l(i)$, $i \in \{0\ldots, J\}$, is the number of the level of the $i$th 3ph-loop function. We note that $\boldsymbol{D}^1_\Lambda$ makes the 3ph-loop functions $H^{1/2}$-stable, while $\boldsymbol{D}^0_\Lambda$ results in $L^2$-stable and $\boldsymbol{D}^2_\Lambda$ in $H^1$-stable 3ph-loop functions. The h-nodal functions, on the other hand, are not $H^{1/2}$-stable [9], and thus, the h-loop functions fail to precondition $\mathcal{T}_\mathrm{A}$. Likewise, the nodal functions are not $H^{1/2}$-stable resulting in $\operatorname{cond}\left(\boldsymbol{\Lambda}^\mathsf{T}\boldsymbol{Z}_\mathrm{A}\boldsymbol{\Lambda}\right) \lesssim 1/h$ [7]. We note that an application of a Jacobi preconditioner can improve the behavior quantitatively, but not qualitatively (see Fig. 1a).

Let $\boldsymbol{\Sigma}_{\mathrm{H}i}$ be the h-n-sol functions for which we define the transformation matrix $\boldsymbol{\Sigma}_\mathrm{H}$. Since $\mathcal{V}$ is inducing an $H^{-1/2}$ inner product, we need a rescaling defined by $\left[\boldsymbol{D}^s_\Sigma\right]_{ii} = 2^{-l(i)}$ to make $\operatorname{div}\boldsymbol{\Sigma}_\mathrm{H}$ $H^{-1/2}$-stable. Then it can be proved that $\operatorname{cond}\left(\boldsymbol{D}_\Sigma\boldsymbol{\Sigma}^\mathsf{T}_\mathrm{H}\boldsymbol{Z}_\Phi\boldsymbol{\Sigma}_\mathrm{H}\boldsymbol{D}_\Sigma\right) \lesssim \log^2\left(1/h^2\right)$ (see [7] and references therein).

We define the overall transformation matrix $\boldsymbol{H}^s_{\boldsymbol{X}/k,\boldsymbol{\Sigma}_\mathrm{H}k} := \left[\boldsymbol{X}\boldsymbol{D}^s_\Lambda/\sqrt{k} \quad \boldsymbol{\Sigma}_\mathrm{H}\boldsymbol{D}_\Sigma\sqrt{k}\right]$, where $\boldsymbol{X}$ can be $\boldsymbol{\Lambda}$, $\boldsymbol{\Lambda}_\mathrm{H}$, or $\boldsymbol{\Lambda}_\mathrm{T}$. Summarizing, we find for the preconditioned EFIE that the loop/h-n-sol basis yields $\operatorname{cond}\left(\left(\boldsymbol{H}^0_{\boldsymbol{\Lambda}/k,\boldsymbol{\Sigma}_\mathrm{H}k}\right)^\mathsf{T}\boldsymbol{Z}_\mathrm{E}\boldsymbol{H}^0_{\boldsymbol{\Lambda}/k,\boldsymbol{\Sigma}_\mathrm{H}k}\right) \lesssim 1/h$ and that the 3ph-loops/h-n-sol basis yields $\operatorname{cond}\left(\left(\boldsymbol{H}^1_{\boldsymbol{\Lambda}_\mathrm{T}/k,\boldsymbol{\Sigma}_\mathrm{H}k}\right)^\mathsf{T}\boldsymbol{Z}_\mathrm{E}\boldsymbol{H}^1_{\boldsymbol{\Lambda}_\mathrm{T}/k,\boldsymbol{\Sigma}_\mathrm{H}k}\right) \lesssim \log^2\left(1/h^2\right)$.

Following these considerations, we study the CFIE. For the sake of brevity, the analysis is carried out only for the conformingly discretized CFIE. The findings are, however, the same as for the standard CFIE. Because of the identity operator of the MFIE, the conditioning of the CFIE is better than that of the EFIE with $\operatorname{cond}\left(\widetilde{\boldsymbol{Z}}_\mathrm{C}\right) \lesssim 1/h$. While the largest singular value of $\mathcal{T}_\Phi$ still grows to infinity, the singular values of $\mathcal{T}_\mathrm{A}$ are shifted from 0 to 1/2 by $\mathcal{I}/2$. Next, we have to study a Helmholtz-decomposed CFIE. Following the argumentation in [16], the compact operator $\mathcal{K}$ can be neglected and it suffices to analyze how $\boldsymbol{G}_\mathrm{ff}$ changes the EFIE behavior.

We notice that when $\mathcal{I}$ is discretized with loop, h-loop and 3ph-loop functions, the resulting matrix is equivalent to the discretization of the Laplace operator in its weak formulation with nodal, h-nodal and 3ph-nodal functions, since $\left(\nabla \times \hat{\boldsymbol{n}}f_i, \nabla \times \hat{\boldsymbol{n}}f_j\right)_{L^2} = \left(\nabla f_i, \nabla f_j\right)_{L^2}$ [16]. The Laplace operator induces an inner product on $H^1/\mathbb{R}$, and hence, we need $H^1$-stable basis functions.

It is well-known that the nodal functions are not $H^1$-stable; the condition number grows with $O(1/h^2)$, and given that the Laplace operator $\Delta_\Gamma$ is a pseudo-differential operator of order +2 and $\mathcal{W}$ is a pseudo-differential operator of order +1, the total order of $\mathcal{W}+\Delta_\Gamma$ is +2, and hence, we have $\operatorname{cond}\left(\boldsymbol{\Lambda}^\mathsf{T}\widetilde{\boldsymbol{Z}}_\mathrm{C}\boldsymbol{\Lambda}\right) \lesssim 1/h^2$ [16]. In other words, loop functions applied to the CFIE result in a conditioning worse than when they are applied to the EFIE. The h-nodal functions are not $H^1$-stable. Yet, for $H^1$ they result in a condition number that grows with $\log^2\left(1/h^2\right)$ [13], and thus, we have $\operatorname{cond}\left(\mathbf{D}_\Lambda^2\mathbf{\Lambda}_\mathrm{H}^\mathrm{T}\widetilde{\mathbf{Z}}_\mathrm{C}\mathbf{\Lambda}_\mathrm{H}\mathbf{D}_\Lambda^2\right) \lesssim 1/h^2$. The 3ph-nodal functions are $H^1$-stable [14], and thus, we have $\operatorname{cond}\left(\mathbf{D}_\Lambda^2\mathbf{\Lambda}_\mathrm{T}^\mathrm{T}\widetilde{\mathbf{Z}}_\mathrm{C}\mathbf{\Lambda}_\mathrm{T}\mathbf{D}_\Lambda^2\right) \asymp 1$. Notice that the matrix $\mathbf{D}_\Lambda^2\mathbf{\Lambda}_\mathrm{T}^\mathrm{T}\mathbf{Z}_\mathrm{E}\mathbf{\Lambda}_\mathrm{T}\mathbf{D}_\Lambda^2$ is ill-conditioned; since we use $H^1$-stable functions, the singular values are now accumulating at zero. This does not destroy the well-conditioning, as $\mathbf{D}_\Lambda^2\mathbf{\Lambda}_\mathrm{T}^\mathrm{T}\mathbf{G}_\mathrm{ff}\mathbf{\Lambda}_\mathrm{T}\mathbf{D}_\Lambda^2$ is spectrally equivalent to $\mathbf{I}$ and thus the spectrum is bounded from below.

As the basis spanned by $\operatorname{div}\boldsymbol{\Sigma}_{\mathrm{H}i}$ has derivative strength when rescaled with $\mathbf{D}_\Sigma$ (since the basis is $H^{-1/2}$-stable), the basis spanned by $\boldsymbol{\Sigma}_{\mathrm{H}i}$ has integrative strength. Given that $\mathcal{I}$ has no derivative terms, the matrix $\mathbf{D}_\Sigma\mathbf{\Sigma}_\mathrm{H}^\mathrm{T}\mathbf{G}_\mathrm{ff}\mathbf{\Sigma}\mathbf{D}_\Sigma$ has singular values clustering around zero, and since $\mathbf{D}_\Sigma\mathbf{\Sigma}_\mathrm{H}^\mathrm{T}\mathbf{Z}_\mathrm{E}\mathbf{\Sigma}\mathbf{D}_\Sigma$ is spectrally equivalent to $\mathcal{I}$, we can conclude that $\operatorname{cond}\left(\mathbf{D}_\Sigma\mathbf{\Sigma}_\mathrm{H}^\mathrm{T}\widetilde{\mathbf{Z}}_\mathrm{C}\mathbf{\Sigma}\mathbf{D}_\Sigma\right) \lesssim \log^2\left(1/h^2\right)$.

Summarizing, we have for the loop/h-n-sol basis preconditioner $\operatorname{cond}\left(\left(\mathbf{H}^0_{\Lambda,\Sigma_\mathrm{H}k}\right)^\mathrm{T}\widetilde{\mathbf{Z}}_\mathrm{C}\mathbf{H}^0_{\Lambda,\Sigma_\mathrm{H}k}\right) \lesssim 1/h^2$, for the h-loop/h-n-sol basis preconditioner $\operatorname{cond}\left(\left(\mathbf{H}^2_{\Lambda_\mathrm{H},\Sigma_\mathrm{H}k}\right)^\mathrm{T}\widetilde{\mathbf{Z}}_\mathrm{C}\mathbf{H}^2_{\Lambda_\mathrm{H},\Sigma_\mathrm{H}k}\right) \lesssim \log^2\left(1/h^2\right)$, and for the 3ph-loop/h-n-sol basis preconditioner $\operatorname{cond}\left(\left(\mathbf{H}^2_{\Lambda_\mathrm{T},\Sigma_\mathrm{H}k}\right)^\mathrm{T}\widetilde{\mathbf{Z}}_\mathrm{C}\mathbf{H}^2_{\Lambda_\mathrm{T},\Sigma_\mathrm{H}k}\right) \lesssim \log^2\left(1/h^2\right)$. Equivalently said, the combination loop/h-n-sol functions does not precondition the CFIE, while both combinations h-loop/h-n-sol and 3ph-loop/h-n-sol basis are a valid preconditioner for the CFIE.

In practical scenarios, the mesh is typically unstructured and thus h-loop and 3ph-loop functions are not available; yet, from the presented theory it is clear that the use of loop functions is not effective and a different strategy is necessary. Define the transformation matrix $\hat{\mathbf{\Lambda}} := \mathbf{\Lambda}\left(\mathbf{\Lambda}^\mathrm{T}\mathbf{\Lambda}\right)^{-1/2}$. Using these orthogonalized loop functions, we obtain as hierarchically preconditioned system

$$\left(\mathbf{H}^0_{\hat{\Lambda},\Sigma_\mathrm{H}k}\right)^\mathrm{T}\widetilde{\mathbf{Z}}_\mathrm{C}\mathbf{H}^0_{\hat{\Lambda},\Sigma_\mathrm{H}k}\mathbf{i} = \left(\mathbf{H}^0_{\hat{\Lambda},\Sigma_\mathrm{H}k}\right)^\mathrm{T}\widetilde{\mathbf{v}}_\mathrm{C}\,, \tag{3}$$

with $\operatorname{cond}\left(\left(\mathbf{H}^0_{\hat{\Lambda},\Sigma_\mathrm{H}k}\right)^\mathrm{T}\widetilde{\mathbf{Z}}_\mathrm{C}\mathbf{H}^0_{\hat{\Lambda},\Sigma_\mathrm{H}k}\right) \lesssim \log^2\left(1/h^2\right)$.

While this approach is theoretically sound, it lacks practical applicability due to the presence of square root matrices. We modify the approach by considering a left instead of the split preconditioner in (3). Thereby, we obtain

$$\mathbf{H}^0_{\hat{\Lambda},\Sigma_\mathrm{H}k}\left(\mathbf{H}^0_{\hat{\Lambda},\Sigma_\mathrm{H}k}\right)^\mathrm{T} = \mathbf{\Lambda}\left(\mathbf{\Lambda}^\mathrm{T}\mathbf{\Lambda}\right)^{-1}\mathbf{\Lambda}^\mathrm{T} + \mathbf{S}\mathbf{D}_\Phi^2\mathbf{S}^\mathrm{T}\,. \tag{4}$$

It has been shown that $\mathbf{P}_\Lambda := \mathbf{\Lambda}\left(\mathbf{\Lambda}^\mathrm{T}\mathbf{\Lambda}\right)^{-1}\mathbf{\Lambda}^\mathrm{T}$ is a projector to the solenoidal Helmholtz subspace [17]. Differently from $\left(\mathbf{\Lambda}^\mathrm{T}\mathbf{\Lambda}\right)^{-1/2}$, the inverse of $\mathbf{\Lambda}^\mathrm{T}\mathbf{\Lambda}$ can be rapidly computed iteratively: the matrix $\mathbf{\Lambda}^\mathrm{T}\mathbf{\Lambda}$ is a graph Laplacian, and black-box preconditioners like algebraic multigrid methods can be used to precondition it [16]. The left preconditioner is motivated by the fact that when we have matrices $\mathbf{A},\mathbf{B} \in \mathbb{R}^{n\times n}$, and $\mathbf{A}$ is symmetric and $\mathbf{x}^\mathrm{T}\mathbf{x} \lesssim \mathbf{x}^\mathrm{T}\mathbf{B}^\mathrm{T}\mathbf{A}\mathbf{B}\mathbf{x} \lesssim \mathbf{x}^\mathrm{T}\mathbf{x}$ forall $\mathbf{x} \in \mathbb{R}^n$, we have $\operatorname{cond}\left(\mathbf{B}\mathbf{B}^\mathrm{T}\mathbf{A}\right) \lesssim 1$. In our case $\mathbf{A} = \widetilde{\mathbf{Z}}_\mathrm{C}$ is not symmetric due to the perturbation by the non-static kernel and the operator $\mathcal{K}$. Since the perturbation is only compact, we don't expect a significant impact on the behavior of the Krylov subspace methods.

Usually, best results are obtained when the (hierarchical) functions are rescaled by leveraging on a Jacobi preconditioner, that is, we use $\left[\mathbf{D}_\Sigma\right]_{ii} = 1/\sqrt{\left[\mathbf{\Sigma}_\mathrm{H}^\mathrm{T}\widetilde{\mathbf{Z}}_\mathrm{C}\mathbf{\Sigma}_\mathrm{H}\right]_{ii}}$. For a far comparison of the different bases discussed, solenoidal and non-solenoidal alike, each basis is rescaled by using such a Jacobi preconditioner.

When the conforming CFIE is used, the Gram matrices prohibit to efficiently obtain it. For the preconditioner presented in this work, this problem can be avoided by using $\left[\mathbf{D}_\Sigma\right]_{ii} = 1/\sqrt{\left[\mathbf{\Sigma}_\mathrm{H}^\mathrm{T}\mathbf{Z}_\mathrm{E}\mathbf{\Sigma}_\mathrm{H}\right]_{ii}}$, that is, the same procedure as for the EFIE can be used. This choice cures the low-frequency breakdown (i.e., the condition is independent of the frequency); however, a further aligning of the singular values branches associated with the solenoidal and the non-solenoidal Helmholtz subspace improves the condition number. To this end, by selecting $\alpha = \left\|\mathbf{P}_\Lambda\widetilde{\mathbf{Z}}_\mathrm{C}\right\|_2$ and $\beta = \left\|\mathbf{S}\mathbf{D}_\Phi\mathbf{D}_\Phi\mathbf{S}^\mathrm{T}\widetilde{\mathbf{Z}}_\mathrm{C}\right\|_2$, we can define $\mathbf{M} := \mathbf{P}_\Lambda/\alpha + \mathbf{S}\mathbf{D}_\Phi\mathbf{D}_\Phi\mathbf{S}^\mathrm{T}/\beta$ resulting in the system

$$\mathbf{M}\widetilde{\mathbf{Z}}_\mathrm{C}\mathbf{i}_\mathrm{M} = \mathbf{M}\widetilde{\mathbf{v}}_\mathrm{C}\,. \tag{5}$$

The norms can be estimated rapidly by using a power iteration method.

## IV. Numerical Results

In the following, we compare our new formulation (5) with the unpreconditioned CFIE and when standard hierarchical basis preconditioners are applied. The standard hierarchical preconditioners we use always consist of the same h-n-sol basis complemented by loop, h-loop, or 3ph-loop functions, respectively. We do this for the conformingly discretized CFIE operator $\widetilde{\mathbf{Z}}_\mathrm{C}$, and the standard CFIE $\mathbf{Z}_\mathrm{C}$ (denoted as S-CFIE).

First, we analyzed the dense-discretization stability by using a cube with side length 1 m. We used a plane wave excitation and fixed the frequency at $f = 1\,\mathrm{MHz}$ with $f = kc/(2\pi)$, where $c$ is the speed of light. We varied the average edge length $h$ from 1.13 m to 0.07 m, and the number of iteration is displayed as a function of the spectral index $1/h$ in Fig. 1a. The results confirm the presented theory: in particular, we can conclude from the figure that loop functions should not be used with the CFIE, and that h-loops can safely be used with the CFIE but not with the EFIE. The new formulation (5) performs well both when applied to a conformingly discretized CFIE as well as a standard CFIE.

Next, we verified the frequency stability. Figure 1b displays the condition number as a function of the frequency $f$. We see that all preconditioned formulations are free from the low-frequency breakdown. Figure 1c shows that all preconditioned CFIEs are resonance-free. We also see that in terms of the condition number, the unpreconditioned CFIE works better than the preconditioned counterparts. This is not unexpected, since hierarchical basis preconditioners need to be adapted to high-frequency problems.

To apply our new method to a more realistic structure, we employed the model of a space shuttle shown in Fig. 2. The electric size of the shuttle is $1/2\lambda$, where $\lambda$ is the wavelength. Since the mesh is unstructured, we do not have a h-sol basis, and thus, we can complement the h-n-sol basis only with loop functions. Table I summarizes our results, where the solver tolerance was $1\times10^{-6}$.

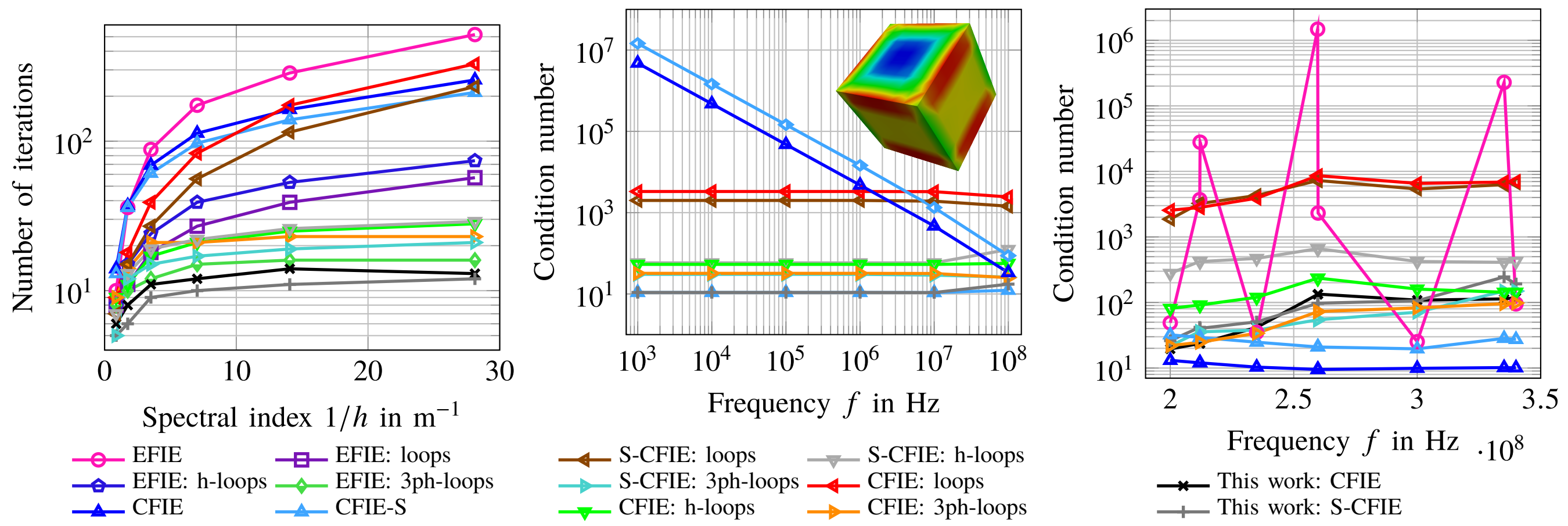

(a) Dense-discretization stability.

(b) Low-frequency stability.

(c) Resonance-freeness.

Fig. 1. Cube: spectral analysis.

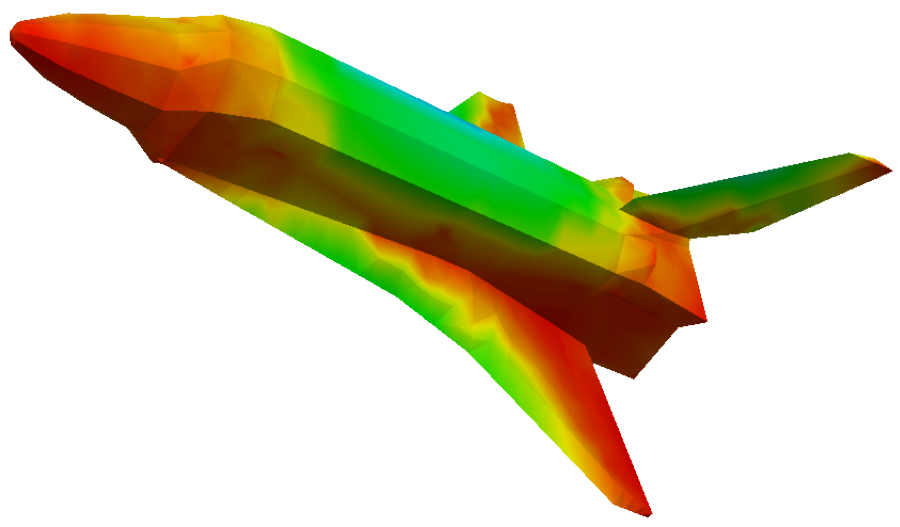

Fig. 2. Shuttle

| Formulation | Iterations |
|---|---|
| CFIE | 400 |
| S-CFIE | 417 |
| CFIE: loops | 253 |
| S-CFIE: loops | 198 |
| This work: CFIE | 33 |
| This work: S-CFIE | 37 |

TABLE I
SHUTTLE: THE NUMBER OF ITERATIONS FOR THE DIFFERENT FORMULATIONS WITH SOLVER TOLERANCE $1 \times 10^{-6}$.

## V. CONCLUSION

First, we can conclude that the h-loop functions, which fail to precondition the vector potential $\mathcal{T}_\mathrm{A}$ of the EFIE, can be successfully applied to the CFIE (though the results obtained by using 3ph-loop functions are usually better, and they work for the EFIE as well). The best results were obtained with our new formulation, which differently from all other hierarchical preconditioners, works on both structured and unstructured meshes.

## REFERENCES

[1] S. H. Christiansen and J.-C. Nédélec, "A Preconditioner for the Electric Field Integral Equation Based on Calderon Formulas," *SIAM J. Numer. Anal.*, vol. 40, no. 3, pp. 1100–1135, Jan. 2002.

[2] W.-L. Wu, A. W. Glisson, and D. Kajfez, "A Study of Two Numerical Solution Procedures for the Electric Field Integral Equation," *Appl. Comput. Electromagn. Soc. J.*, vol. 10, no. 3, pp. 69–80, 1995.

[3] J. R. Mautz and R. F. Harrington, "H-Field, E-Field, and Combined-Field Solutions for Conducting Bodies of Revolution," *Archiv für Elektronik und Übertragungstechnik (A. E. U)*, vol. 32, pp. 157–163, 1978.

[4] F. Vipiana, P. Pirinoli, and G. Vecchi, "A Multiresolution Method of Moments for Triangular Meshes," *IEEE Trans. Antennas Propag.*, vol. 53, no. 7, pp. 2247–2258, Jul. 2005.

[5] F. P. Andriulli, F. Vipiana, and G. Vecchi, "Hierarchical Bases for Non-hierarchic 3-D Triangular Meshes," *IEEE Trans. Antennas Propag.*, vol. 56, no. 8, pp. 2288–2297, Aug. 2008.

[6] R.-S. Chen, J. Ding, D. Ding, Z. Fan, and D. Wang, "A Multiresolution Curvilinear Rao-Wilton-Glisson Basis Function for Fast Analysis of Electromagnetic Scattering," *IEEE Trans. Antennas Propag.*, vol. 57, no. 10, pp. 3179–3188, 2009.

[7] F. P. Andriulli, A. Tabacco, and G. Vecchi, "Solving the EFIE at Low Frequencies with a Conditioning That Grows Only Logarithmically with the Number of Unknowns," *IEEE Trans. Antennas Propag.*, vol. 58, no. 5, pp. 1614–1624, May 2010.

[8] M. A. Francavilla, M. Righero, F. Vipiana, and G. Vecchi, "Applications of a Hierarchical Multiresolution Preconditioner to the CFIE," in *Eur. Conf. on Antennas and Propagation*, The Hague, Netherlands, Apr. 2014.

[9] R. Lorentz and P. Oswald, "Multilevel finite element Riesz bases in Sobolev spaces," *DD9 Proceedings, Domain Decomposition Press, Bergen*, pp. 178–187, 1998.

[10] M. Righero, I. M. Bulai, M. A. Francavilla, F. Vipiana, M. Bercigli, A. Mori, M. Bandinelli, and G. Vecchi, "Hierarchical Bases Preconditioner to Enhance Convergence of the CFIE with Multiscale Meshes," *Posted on arXiv*.

[11] S. B. Adrian, F. P. Andriulli, and T. F. Eibert, "Hierarchical Bases Preconditioners for a Conformingly Discretized Combined Field Integral Equation Operator," in *Eur. Conf. on Antennas and Propagation*, Lisbon, Portugal, Apr. 2015.

[12] K. Cools, F. P. Andriulli, D. De Zutter, and E. Michielssen, "Accurate and Conforming Mixed Discretization of the MFIE," *IEEE Antennas Wireless Propag. Lett.*, vol. 10, pp. 528–531, 2011.

[13] H. Yserentant, "On the Multi-Level Splitting of Finite Element Spaces," *Numer. Math.*, vol. 49, no. 4, pp. 379–412, Jul. 1986.

[14] R. Stevenson, "Stable Three-Point Wavelet Bases on General Meshes," *Numer. Math.*, vol. 80, no. 1, pp. 131–158, Jul. 1998, 00029.

[15] J.-C. Nédélec, *Acoustic and Electromagnetic Equations: Integral Representations for Harmonic Problems*, ser. Applied Mathematical Sciences. New York: Springer, 2001, no. 144.

[16] F. Andriulli, "Loop-Star and Loop-Tree Decompositions: Analysis and Efficient Algorithms," *IEEE Trans. Antennas Propag.*, vol. 60, no. 5, pp. 2347–2356, May 2012.

[17] F. P. Andriulli, K. Cools, I. Bogaert, and E. Michielssen, "On a Well-Conditioned Electric Field Integral Operator for Multiply Connected Geometries," *IEEE Trans. Antennas Propag.*, vol. 61, no. 4, pp. 2077–2087, Apr. 2013.